\documentclass{elsart}

\usepackage{amssymb}
\usepackage{amsmath}

\numberwithin{equation}{section} \numberwithin{thm}{section}

\begin{document}

\begin{frontmatter}


 \title{Ramanujan and Eckford Cohen totients from Visible Point Identities}
\thanks[label1]{Thanks are due to the La Trobe University Mathematics Department for their affiliation during some of the research for this paper.}
 \author{Geoffrey B Campbell}
 \ead{Geoffrey.Campbell@anu.edu.au}
 \address{Mathematical Sciences Institute \\ The Australian National University \\ ACT, 0200, Australia}





\begin{abstract}
We define an extension of Ramanujan's trigonometric function to
arbitrary dimensions, and give the Dirichlet series generating
function.  The extension was first given by Eckford Cohen long ago. 
This links directly to visible point vector (vpv)
identities in the author's papers, and possibly to lattice sums in
Physics and Chemistry presented by Baake \emph{et al}. New
generating functions and summations are given here, generalizing
the Ramanujan function, Euler totient and the Jordan totient
functions, based on visible lattice point ideas.

\end{abstract}


\end{frontmatter}

\section{Introduction}
\label{s:1}

    The function made well known by Ramanujan~\cite{sR1962} for positive integers $k$, and~$n$,

\begin{equation}
c_{k}(n)
= \sum_{ \substack{(j,k)=1 \\ 0<j<k}} \cos \left(\frac{2 \pi
nj}{k}\right) =\sum_{\substack{(j,k)=1\\0<j<k}} e^{\frac{2\pi
inj}{k}}
\end{equation}

has been widely studied and referenced.  See for example
Hardy~\cite{gH1940}, Hardy and Wright~\cite{gH1971},
Sivaramakrishnan~\cite{rS1989}, and Apostol~\cite{tA1976}.  The paper by Cohen~\cite{eC1959} contains versions of some of the formulas of the current paper, couched in the language of arithmetical functions of the 1950s. The perspective offered in the current paper shows a natural approach to these identities, with its genesis in the visible point vector identities.  In this note we extend the function as follows

\begin{equation}
c_{k}(n_1,n_2)
= \sum_{ \substack{(j_1,j_2,k)=1 \\ 0\leq j_1<k\\0\leq j_2<k\\
j_1+j_2\neq 0}} \cos \left(\frac{2 \pi(j_1n_1+j_2n_2)}{k}\right)
=\sum_{\substack{(j_1,j_2,k)=1 \\ 0\leq j_1<k\\0\leq j_2<k\\
j_1+j_2\neq 0}} \exp\left(\frac{2\pi i(j_1n_1+j_2n_2}{k}\right),
\end{equation}

\begin{multline}
c_{k}(n_1,n_2, n_3)
= \sum_{ \substack{(j_1,j_2, j_3, k)=1 \\ 0\leq j_1<k\\0\leq j_2<k\\
0\leq j_3<k\\j_1+j_2+j_3\neq 0}} \cos \left(\frac{2
\pi(j_1n_1+j_2n_2+j_3n_3)}{k}\right)\\
=\sum_{\substack{(j_1,j_2,j_3,k)=1 \\ 0\leq j_1<k\\0\leq j_2<k\\
0\leq j_3<k \\ j_1+j_2+j_3\neq 0}} \exp\left(\frac{2\pi
i(j_1n_1+j_2n_2+j_3n_3}{k}\right),
\end{multline}

and generally

 \begin{multline}
c_{k}(n_1,n_2, \dots , n_m)
= \sum_{ \substack{(j_1,j_2, \dots, j_m, k)=1 \\ 0\leq j_1<k\\\vdots\\0\leq j_m<k\\
j_1+j_2+\dots +j_m\neq 0}} \cos \left(\frac{2
\pi(j_1n_1+j_2n_2+\dots +j_mn_m)}{k}\right)\\
=\sum_{\substack{(j_1,j_2, \dots, j_m,k)=1 \\ 0\leq j_1<k\\ \vdots\\ 0\leq j_m<k\\
j_1+j_2+ \dots +j_m\neq 0}} \exp\left(\frac{2\pi
i(j_1n_1+j_2n_2+\dots +j_mn_m}{k}\right)
\end{multline}

These functions are very natural to consider as their generating
Dirichlet series exhibit a simplicity unexpected of a function
requiring so many parameters. On the other hand, there seems an
inherently natural aspect to considering the set of visible
lattice points in any radial region of multidimensional Euclidean
space.  We shall see in the following sections that
$c_{k}(n_1,n_2, \dots , n_m)$ always reduces easily to a one
dimensional  $c_k(n)$   type function.

\section{Dirichlet series generating functions}

In this section we start with a definition then give the
generating function for the new Ramanujan function, and some
consequences.

    \begin{defn}  Let us define the operator

\begin{equation}
\sum_{\substack{(j_1j_2, \dots,j_m,k)=1\\ 0\leq j_1<k\\\vdots\\
0\leq j_m<k\\ j_1+j_2+\dots +j_m \neq 0}} \substack{:}{=}
{\sum_m}_k,
\end{equation}

so that

 \begin{equation}
 c_k(n_1, n_2, \dots, n_m)={\sum_{m}}_k
 \cos\left(\frac{2\pi(j_1n_1+j_2n_2+\dots +j_mn_m)}{k}\right).
 \end{equation}

    \end{defn}

    It is well known that~\cite[pp. 139-143]{gH1971}

\begin{equation}
\frac{\sigma_{-s}(n)}{\zeta(s+1)} =\sum^\infty_{k=1}
\frac{c_k(n)}{k^{s+1}}, \qquad \mathcal{R}s>0,
\end{equation}

where the left side contains the sum of -$s$th  powers of the
divisors of $n$, and the Riemann zeta function.  If the Riemann
hypothesis is true then (2.3) is valid for
$\mathcal{R}s>-\frac{1}{2}$.  We show that
    \begin{thm}  If  m  is a positive integer then,

\begin{equation}
\frac{\sigma_{-s-1+m}((n_1, n_2, \dots, n_m))}{\zeta(s+1)}
=\sum^\infty_{k=1} \frac{c_k(n_1, n_2, \dots, n_m)}{k^{s+1}},
\qquad \mathcal{R} s>0,
\end{equation}

 \begin{equation}
\exp\left(\sum_{k\mid(n_1, \dots, n_m)} k^{m-1}z^k\right)
=\prod^\infty_{k=1} (1-z^k)^{-c_{k}(n_1, \dots, n_m)/k}, \qquad
\vert z\vert <1.
 \end{equation}

 \end{thm}

Like (2.3), (2.4) is probably also true for
$\mathcal{R}s>-\frac{1}{2}$.  It is easy to prove many results
similar to (2.4) and (2.5).  We shall prove theorem 2.2 in the
next section.  Features of the Ramanujan function such as
multiplicativity generalize in the new versions.
    \begin{cor}  $c_k(n_1, n_2, \dots, n_m)$  is  multiplicative in  k. \end{cor}
The paper by Apostol and Zuckermann~\cite{tA1964} deals with the
functional equation  $F(mn) F((m, n)) = F(m) F(n) f((m,n))$  which
many of the common arithmetical functions satisfy, and indeed so
do our new functions   $c_k(n_1, n_2, \dots, n_m)$, in the form
\begin{multline}
c_{mn}(n_1, n_2, \dots, n_\lambda)c_{(m,n)}(n_1, n_2, \dots,
n_\lambda)\\ =c_{m}(n_1, n_2, \dots, n_\lambda) \sum^\infty_{k=1}
\frac{c_k(n_1, n_2, \dots, n_m)}{k} =0 f((m,n)).
 \end{multline}

    \begin{cor}  If  $m$  is a positive integer,

\begin{equation}
\sum^\infty_{k=1} \frac{c_k(n_1, n_2, \dots, n_m)}{k} =0.
\end{equation}

    \end{cor}

Following the original paper by Ramanujan~\cite{sR1962} on
$c_{n}(k)$  this comes from the fact that the ``\emph{sum of
powers of divisors}'' function on left of (2.4) is a finite
Dirichlet series, and at  $s=0$  in that equation depends on the
well known convergence of

\begin{equation}
\sum^\infty_{k=1} \frac{\mu(k)}{k}=0.
\end{equation}

\section{Proof of Theorem 2.2}

Theorem 2.2 depends on lemmas 3.1 and 3.2.  The latter is an
analytic restatement of the former.
\begin{lem} (see Campbell~\cite{gC1994a})  Consider an infinite
region raying out of the origin in any Euclidean vector space. The
set of all lattice point vectors from the origin in that region is
precisely the set of positive integer multiples of the visible
point vectors  (vpv's) in that region.
\end{lem}
This lemma underlies all of the author's cited
papers~\cite{gC1998} to \cite{gC1998a} and is fundamental in
obtaining vpv identities.  A clearer picture of lemma 3.1 is
gained from fig 1.  It shows part of the 2-D first quadrant
visible-from-origin lattice points as the dots.  All other lattice
points in the part quadrant are shown as crosses.  By lemma 3.1,
the coordinates of the crosses are the positive integer (scalar)
multiples of the coordinates of the dots; this being true in fact
for the dots and crosses enclosed in any radial region from the
origin. In the study of quasicrystals Baake~et~al~\cite{mB1994} for example, diffraction patterns and sectioning gives rise to Fourier
inversions of lattice patterns such as the visible points depicted
here. The identities in section 3 here are an "easy" way into many 
arithmetical sums not unlike those from Ramanujan~\cite{sR1962}.

\begin{center}
\begin{tabular}{c|ccccccccc}
  $\boldsymbol {Y}$ &&&&&&&&\\
  8 \ & $\bullet$ & x & $\bullet$ & x & $\bullet$ & x & $\bullet$ & x &\\
  7 \ & $\bullet$ & $\bullet$& $\bullet$ & $\bullet$ & $\bullet$ & $\bullet$ & x & $\bullet$ &\\
  6 \ & $\bullet$ & x & x & x & $\bullet$ & x & $\bullet$ & x &\\
  5 \ & $\bullet$ & $\bullet$ & $\bullet$ & $\bullet$ & x & $\bullet$& $\bullet$ & $\bullet$ &\\
  4 \ & $\bullet$ & x &   $\bullet$ & x& $\bullet$ & x & $\bullet$ &  x&\\
  3 \ & $\bullet$ & $\bullet$ & x & $\bullet$ & $\bullet$ & x & $\bullet$ & $\bullet$ &\\
  2 \ & $\bullet$ & x & $\bullet$ & x & $\bullet$ & x & $\bullet$ & x &\\
  1 \ & $\bullet$ & $\bullet$ & $\bullet$ & $\bullet$ & $\bullet$ & $\bullet$ & $\bullet$ & $\bullet$ &\\ \hline
  0 \ &1&2&3&4&5&6&7&8& $\boldsymbol{X}$
\end{tabular}

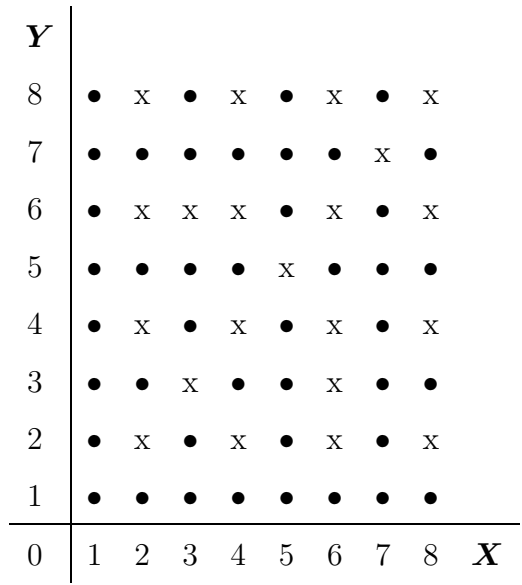
fig. 1.\end{center}


An analytic restatement of lemma 3.1 is given in the
\begin{lem}   If  $(a_k)$ is an arbitrary sequence and  $q_i$  are
variables chosen so that the following functions are all defined,

\begin{equation}
\sum^\infty_{k=1} \left(a_k \prod^m_{h=1}
\frac{1-q_h}{1-q_h^{1/k}}\right)=S_1+\sum^\infty_{k=2}
\left(S_k{\sum_{m}}_k (q_1^{j_{1}} \dots q^{j_{m}}_m )^{1/k}
\right) \text{ where } S_k =\sum^\infty_{j=1} a_{jk}.
\end{equation}

\end{lem}

\textbf{Proof of Lemma 3.2.}   The case with  $m=1$  was given in
Campbell~\cite{gC1998} to~\cite{gC1998a}, leading
in~\cite{gC1992} to proofs of (2.3) and the case  $m=1$  of
(2.5).  The cases of (3.1) with   $m = 1, 2, 3,$ are as follows,
and correspond to lemma 3.1 in $2$-D, $3$-D, and $4$-D space :-

\begin{equation}
\sum^\infty_{k=1} a_k \frac{1-q_1}{1-q_1^{1/k}}
=S_1+\sum^\infty_{k=2} \left(S_k {\sum_1}_k q_1^{j_{1}/k}\right),
\end{equation}

 \begin{equation}
\sum^\infty_{k=1} a_k \frac{1-q_1}{1-q_1^{1/k}}
\frac{1-q_2}{1-q_2^{1/k}} =S_1+\sum^\infty_{k=2} \left(S_k
{\sum_{2}}_k (q_1^{j_{1}}q_2^{j_{2}})^{1/k}\right),
\end{equation}

\begin{equation}
\sum^\infty_{k=1} a_k \frac{1-q_1}{1-q_1^{1/k}}
\frac{1-q_2}{1-q_2^{1/k}}\frac{1-q_3}{1-q_3^{1/k}}
=S_1+\sum^\infty_{k=2} \left(S_k
{\sum_{3}}_k(q_1^{j_{1}}q_2^{j_{2}}q_3^{j_{3}})^{1/k}\right),
\end{equation}

The proof of each of (3.2) to (3.4) is almost \emph{a priori} when
seen as interpreting lemma 3.1 as it is depicted in a generalized
version of fig 1.  To spell this out, we see that the left side of
(3.1) is

\begin{align*}
&\quad  \sum^\infty_{k=1} \left(a_k\prod^m_{h=1}
\frac{1-q_h}{1-q_h^{1/k}}\right)\\
&=\sum^\infty_{k=1} \left(a_k\prod^m_{h=1} \sum^{k-1}_{j=0}
q^{j/k}_h\right)\\
&=\sum^\infty_{k=1} \left(a_k\sum^{k-1}_{j=1}(q_1^{j_{1}} \dots
q_m^{1/k})\right)\\
&=S_1+\sum^\infty_{k=2} \left(S_k{\sum_{m}}_k(q_1^{j_{1}} \dots
q_m^{1/k})\right) \text{ where } S_k =\sum^\infty_{j=1} a_{jk}.
\end{align*}

The last step here depends on application of lemma 3.1 in
$m$-space. $\blacksquare$

The vpv identities (visible point vector) published by the author,
such as for example in Campbell~\cite{gC1994a} or~\cite{gC2000}
are all particular cases of lemma 2.2.

\textbf{Proof of Theorem 2.2.}  In lemma 3.2  set  $q_k=exp(2\pi i
n_k)$, then set $a_k$  equal to respectively  $k^{-s}$,  $z^k$.
The resulting identities are (2.4) and (2.5).
    $\blacksquare$

\section{A new Jordan Totient generating function, and some related
results}

The consequent identities from (3.2) were explored in some detail
in the author's papers~\cite{gC1998} to~\cite{gC1998a}.  Borwein
and Borwein~\cite[ex 10, pp.327]{jB1986} had some infinite
products resembling those of the author, and some of the Dirichlet
series implied from (3.2) were given in Lossers~\cite{oL1985},
Sivaramakrishnan~\cite{rS1989}, and
Chandrasekharan~\cite{kC1970}.  Most books on arithmetical
functions since Ramanujan's original paper~\cite{sR1962} have
included a treatment of the function (1.1) and give at least the
result (2.3) cited above.  However, (2.4) and (2.5) are not in the
literature, and the identities (3.3) and (3.4) for example, lead
us to a class of Dirichlet series and infinite products not
hitherto considered.  These results are near the surface and seem
worth closer study.  Lemma 3.2 is based on a summation over
lattice points in an  $m+1$  dimensional hyperpyramid.  (see
Campbell~\cite{gC2000})  The methods for lattice sums given in
Ninham~et~al~\cite{bN1992} and Glasser and Zucker~\cite{mG1980}
using Mobius and Mellin inversions are also applicable to results
from this lemma. Of course this has been recently given a more
systematic treatment by Baake~et~al~\cite{mB1994}
to~\cite{mB2002a} and his colleagues in their investigations of
quasicrystal lattice structures and associated Fourier inversions.

Since some neat examples of (3.2) have been examined in other
papers, we now consider the next simplest equation (3.3) with
$q_1, q_2$ replaced respectively by $e^{xz}$, $e^{yz}$. We have
therefore the analysis

 \begin{equation}
\sum^\infty_{k=1} a_k \left(\sum^{k-1}_{A=0} \sum^{k-1}_{B=0}
e^{(Ax+By)z/k}\right) =S_1 + \sum^\infty_{k=2}
\left(S_k{\sum_{2}}_k e^{(j_1x+j_2y)z/k}\right),
 \end{equation}

 \begin{equation}
\Leftrightarrow \sum^\infty_{k=1} a_k \left(\sum^{k-1}_{A=0}
\sum^{k-1}_{B=0} \sum^{k-1}_{C=0}\frac{(Ax+By)^c z^c}{C!
k^c}\right) =S_1 + \sum^\infty_{k=2} \left(S_k{\sum_{2}}_k
 \sum^\infty_{C=0} \frac{(j_1x+j_2y)^cz^c}{C!
k^c}\right),
 \end{equation}

\begin{equation}
\Leftrightarrow \left(\sum^{\infty}_{k=1}k^2a_k\right)
+\sum^{\infty}_{C=0} \sum^{\infty}_{k=1}\sum_{\substack{0\leq
A<k\\ 0\leq B<k\\A+B\neq 0}}\frac{(Ax+By)^c a_k z^c}{k^c C!}
 \end{equation}
$$=S_1+\sum^\infty_{C=0} \sum^\infty_{k=2}\left(S_k {\sum_{2}}_k
\frac{(j_1 x+ j_2 y)^c}{k^c} \right) \frac{z^c}{C!},$$

 where, equating coefficients of like powers of $z$ gives us the summation
formulae

\begin{equation}
\sum^\infty_{k=1} k^2 a_k=S_1+\sum^\infty_{k=2} S_k {\sum_{2}}_k
1,
\end{equation}

 \begin{equation}
\sum^\infty_{k=1} \sum_{\substack{0\leq A<k \\0\leq B<k\\A+B
\neq0}} \frac{(Ax+By)a_k}{k} =\sum^\infty_{k=2}
\left(\frac{S_k}{k} {\sum_{2}}_k (j_1x+j_2 y)\right),
 \end{equation}

 \begin{equation}
\sum^\infty_{k=1} \sum_{\substack{0\leq A<k \\0\leq B<k\\A+B
\neq0}} \frac{(Ax+By)^2 a_k}{k^2} =\sum^\infty_{k=2}
\left(\frac{S_k}{k^2} {\sum_{2}}_k (j_1x+j_2 y)^2\right),
 \end{equation}

and so on.  If this process from (4.1) to (4.6) is followed in an
analogous manner for the general series (3.1), we have
\begin{cor}
For positive integers $n$,

 \begin{equation}
\sum^\infty_{k=1} \sum_{\substack{0\leq A<k \\0\leq B<k\\A+B
\neq0}} \frac{(Ax+By)^n a_k}{k^n} =\sum^\infty_{k=2}
\left(\frac{S_k}{k^n} {\sum_{2}}_k (j_1x+j_2 y)^n\right).
 \end{equation}
\end{cor}
In (4.4) the term  $${\sum_{2}}_k 1={\sum_{2}}_k
\left(\frac{j_1+j_2}{k}\right)^0=\varphi_0(2; k)$$   is the number
of non-negative integer solutions of $(a,b,k)=1$  for  $a$ and $b$
both less than  $k$ with $a+b\neq0$. Suggestive from this, the
following is a multidimensional totient function quite distinct
from those of section 1.
\begin{defn}
  For all positive integers  $k$,  $m$,  and  $t + 1$,  let

 \begin{equation}
\varphi_t(m; k)={\sum_{m}}_k \left(\frac{j_1+\dots
+j_m}{k}\right)^t.
 \end{equation}
\end{defn}

If  $t$  is zero and  $m=1$  we would have the Euler totient
function.  If  $t$  is zero and  $m=2$  we have again the function
$\varphi_0(2; k)$. We shall see soon that  $\varphi_0(m; k)$ is
the Jordan totient function $J_m(k)$. The Dirichlet series
generating functions for $\varphi_t(m; k)$ are found in a similar
fashion to our derivation of corollary 4.1, and setting $x=y=1$.
If this analogous process is followed we arrive at the
    \begin{thm}   For positive integers $t$   and suitably chosen functions $a_i$,

 \begin{equation}
\sum^\infty_{k=1} k^ma_k=S_1+\sum^\infty_{k=2} S_k \;
\varphi_t(m;k),
 \end{equation}

whilst for $t=0$ we have

\begin{equation}
\sum^\infty_{k=1} k^m a_k=S_1+\sum^\infty_{k=2} S_k \;
\varphi_0(m;k).
 \end{equation}
\end{thm}

    Therefore if  $a^k=k^{-z}$   we have
\begin{cor}

\begin{equation}
\frac{\zeta (s-m)}{\zeta(s)}=1+\sum^\infty_{k=2}
\frac{\varphi_0(m;k)}{k^s}, \quad \mathcal{R} s>1.
\end{equation}
\end{cor}

So it is clear that  $\varphi_0(m;k)=J_m(k)$, the Jordan totient
function. This function is well known as the number of ordered
sets of $k$  elements chosen from a complete residue system  (mod
r) such that the greatest common divisor of each set is prime to
$r$.  It is also well known (see Sivaramakrishnan~\cite{rS1989}),
and easily deduced from writing left side of (4.11) as a product
over primes, that

\begin{equation}
J_m(k) =k^m\prod_{p\vert k}(1-p^{-m}).
\end{equation}

This function is therefore perhaps more simply defined as the
number of non-negative integer solutions of
$(a_1,a_2,\ldots,a_n,k)=1$ with all of the $a$'s less than $k$,
and $a_1+a_2+\cdots+a_n \neq 0$.
    If in (4.9) we let  $a^k=z^k k^{-1}$ then

\begin{cor}

 \begin{equation}
\frac{1}{1-z} \prod_{k=2}\left(\frac{1}{1-z^k}\right)^{J_{m}(k)/k}
=\exp \left\{\sum^\infty_{k=1} k^{m-1} z^k\right\}, \quad \vert
z\vert <1.
 \end{equation}
\end{cor}

This is new, and related to many results in the author's papers
cited, and also to a product in Borwein and Borwein~\cite[ex 10,
pp.327]{jB1986}.  The right side of (4.13) is easily reduced by
using the well known result

 \begin{equation}
\sum^\infty_{k=0} k^{m-1} z^k=\sum^{m-1}_{j=0} S_{m-1}^{(j)} z^j
\frac{d^j}{dz^j} \left\{ \frac{1-z^n}{1-z}\right\}, \quad z\neq 1,
\end{equation}

where $S_{m-1}^{(j)}$ are the Stirling numbers of the second kind.
Hence we have the
\begin{thm}
If  $m$  is a positive integer

\begin{equation}
\frac{1}{1-z} \prod^\infty_{k=2}
\left(\frac{1}{1-z^k}\right)^{J_{m}(k)/k}=\exp
\left\{\sum^{m-1}_{j=0} S^{(j)}_{m-1}
\frac{j!z^j}{(1-z)^{j+1}}\right\}, \quad \vert z \vert <1.
\end{equation}
\end{thm}

The Stirling numbers of the second kind are of importance in
combinatorial number theory. The Stirling number $S^{(n)}_{m}$ is
the number of ways of partitioning a set of $m$ elements into $n$
non-empty sets. In light of (4.13), formula (4.15) corresponds to
a limiting case of the hyperpyramid lattice vpv identity first
given in Campbell~\cite{gC1994a},

 \begin{equation}
\prod_{\substack{(a_1, a_2, \dots, a_n)=1\\ a_1, \dots,
a_{n-1}<a_n \\
a_1, \dots, a_{n-1} \geq 0\\
a_n \geq 1}} \left(\frac{1}{1-x^{a_{1}}_1 x_2^{a_{2}} \dots
x_n^{a_{n}}}\right)^{\frac{1}{a_1^{b_{1}} a_2^{b_{2}} \dots
a_n^{b_{n}}}} =\exp \left\{ \sum^\infty_{k=1} \prod^{n-1}_{i=1}
\left(\sum^{k-1}_{j=1} \frac{x^j_i}{j^{b_{i}}}\right)
\frac{x^k_n}{k^{b_{n}}}\right\}
 \end{equation}

where for  $i = 1, 2, ..., n$;  $|x_i| < 1$,  $b_i$  is a complex
number with $\sum^n_{i=1} b_i=1$.

\section{Further multidimensional formulae}
    In Campbell [8] we proved the
\begin{thm}
If  $(a_i)$, $(b_i)$   are arbitrary sequences chosen so that,
    together with choice of $x$,  the following functions are defined then
\begin{equation}
\sum^n_{k=1}a_k \frac{1-\exp(b_k x)}{1-\exp (b_k x/k)}
=\left(\sum^n_{k=1} a_k\right)+\sum^n_{m=2}\sum^{[n/m]}_{k=1}
\sum_{\substack{(j,m)=1\\ 0<j<k}} a_{mk} \exp (b_{mk} jx/k),
\end{equation}
where  $[n]$  denotes the greatest integer in  $n$.
\end{thm}

    This theorem has led to many new results on vpv identities, namely:-

a)  the so-called ``\emph{companion identities}'', where a group
of interrelated identities divide up space radially in a
``\emph{raying from origin}'' region.  These were often infinite
products in the examples given.

b)  new results on Dirichlet series.  These included generating
functions for Ramanujan type trigonometric functions $c_k(n)$.

c)  new results involving Jacobi theta functions.

d)  new results on Diophantine equations.

In this section we give further results like (5.1) and some
corresponding results as a consequence.  Theorem 5.1 is a 2-D vpv
summation formula, and it was seen in Campbell~\cite{gC1998}
to~\cite{gC1998a} that vpv identities can arise in any dimension.
This is true in any particular radial region of space with apex at
the origin.  We now give the 3-D version of theorem 5.1.

\begin{thm}
   If   $(a_i)$, $(b_i)$, $(c_i)$,    are each arbitrary sequences chosen so that,
    together with choice of $x$,  the following functions are defined then

\begin{equation} \notag
\sum_{k=1}^n a_k \frac{1-\exp(b_kx)}{1-\exp(b_kx/k)}
\frac{1-\exp(c_kx)}{1-\exp(c_kx/k)}
\end{equation}
\begin{equation} \notag
=\left(\sum^n_{k=1} a_k\right) +\sum^{[n/2]}_{k=1} {\sum_{2}}_2
a_{2k} \exp\left(b_{2k}j_1+c_{2k}j_2)\frac{x}{2}\right)
\end{equation}
\begin{equation}
+\sum^{[n/3]}_{k=1} {\sum_{2}}_3 a_{3k}
\exp\left((b_{3k}j_1+c_{3k}j_2)\frac{x}{3}\right) \end{equation}
\begin{equation} \notag
+\sum^{[n/4]}_{k=1} {\sum_{2}}_4 a_{4k}
\exp\left((b_{4k}j_1+c_{4k}j_2)\frac{x}{4}\right) \end{equation}
\begin{equation} \notag
\vdots \end{equation}
\begin{equation} \notag
+\sum^{[n/n]}_{k=1} {\sum_{2}}_n a_{nk}
\exp\left((b_{nk}j_1+c_{nk}j_2)\frac{x}{n}\right)
\end{equation}

where  $[n]$  denotes the greatest integer in  $n$.
\end{thm}

Most of the 3-D vpv infinite products found in Campbell [10] are
also corollaries of theorem 5.2.  For example, if
$n\rightarrow\infty, a_k=z^k/k, b_k=(k/x)logx, c_k~=~(k/x)~logy$
we have the
\begin{cor} If  $a$  and  $b$  are non-negative
integers and  $c$
     positive integers such that  $(a, b, c) = 1$ then

 \begin{equation}
\prod_{a,b<c}(1-x^ay^b z^c)^{-1/c}
=\left(\frac{(1-xz)(1-yz)}{(1-z)(1-xyz)}\right)^{\frac{1}{(1-x)(1-y)}}
 \end{equation}

valid for every one of $|x|, |xz|, |yz|, |xyz| < 1$.
\end{cor}

This identity was given in [10] as one of four ``companion
identities'' dividing up the first hyperquadrant in 3-D space into
three radial regions and the union region of these.  Of course the
simplicity of statement of (5.3) is slightly disguised in the
notation of theorem 5.2.  Also, the cases of theorem 2.2 with
$m~=~2$  are trivially obtained from theorem 5.2.  Equation 3.2 is
also an easy consequence of (5.2).  If in theorem 5.2 we let both
$b_k$ and  $c_k$  approach unity, we get the new result,

\begin{cor}
If $(a_k)$ is an arbitrary sequence,

 \begin{equation}
\sum^n_{k=1} a_kk^2 = \left(\sum^k_{k=1} a_k\right)+ \left(
\sum^{[n/2]}_{k=1} a_{2k} J_2(2)\right)+ \left( \sum^{[n/3]}_{k=1}
a_{3k} J_2(3)\right)+ \dots + \left( \sum^{[n/n]}_{k=1} a_{nk}
J_2(n)\right).
 \end{equation}
\end{cor}

    An obvious example is

\begin{equation}
\frac{n(n+1)(2n+1)}{6} =n+[n/2]J_2(2) +[n/3]J_2 (3) +\dots +
[n/n]J_2(n).
\end{equation}

It is easy to follow an analogous line of reasoning to that
deriving theorem 5.2, then to obtain the analogue to corollary
5.4, which we rate here as a theorem,
\begin{thm}If  $(a_k)$ is an
arbitrary sequence and  $m$  any positive integer,

 \begin{equation}
\sum^n_{k=1} a_k k^m=\left(\sum^n_{k=1} a_k\right) +J_m(2)
\left(\sum^{[n/2]}_{k=1} a_{2k}\right) +J_m(3)
\left(\sum^{[n/3]}_{k=1} a_{3k}\right) +\dots + J_m(n)
\left(\sum^{[n/n]}_{k=1} a_{nk}\right).
 \end{equation}
\end{thm}

We therefore see that (5.5) is a case of this.  The known
Dirichlet series generating function for the Jordan totient
function is seen to be a limiting case of (5.6) if  $a_k =
k^{-(s+m)}$ and the appropriate convergence restrictions are in
place. Partial sums of this generating function are seen in
\begin{thm} If $m$ and  $n$  are positive integers

\begin{equation}
n=\left(\sum^n_{k=1} \frac{1}{k^{m}}\right) +\frac{J_m(2)}{2^m}
\left(\sum^{[n/2]}_{k=1} \frac{1}{k^{m}}\right) +
\frac{J_m(3)}{3^m} \left(\sum^{[n/3]}_{k=1}
\frac{1}{k^{m}}\right)+\dots + \frac{J_m(n)}{n^m}
\left(\sum^{[n/n]}_{k=1} \frac{1}{k^{m}}\right),
\end{equation}
\begin{multline}
\sum^n_{k=1} k^a=\left(\sum^n_{k=1} \frac{1}{k^{m-a}}\right)
+\frac{J_m(2)}{2^{m-a}} \left(\sum^{[n/2]}_{k=1}
\frac{1}{k^{{m-a}}}\right) \\+ \frac{J_m(3)}{3^{m-a}}
\left(\sum^{[n/3]}_{k=1} \frac{1}{k^{{m-a}}}\right)+\dots +
\frac{J_m(n)}{n^{m-a}} \left(\sum^{[n/n]}_{k=1}
\frac{1}{k^{{m-a}}}\right),
\end{multline}
\begin{equation}
\sum^n_{k=1} k^m=n+[n/2]J_m (2) +[n/3]J_m(3) +\dots + [n/n]J_m(n).
\end{equation}
\end{thm}

\begin{cor}  If  $m$  and  $n$  are positive integers with  $|z| <
1$
 \begin{equation}
\sum^n_{k=1} z^kk^m =\frac{1-z^n}{1-z} +J_m(2)
\frac{1-z^{2[n/2]}}{1-z^2}+J_m(3) \frac{1-z^{3[n/3]}}{1-z^3}
+\dots + J_m(n) \frac{1-z^{n[n/n]}}{1-z^n}.
 \end{equation}
\end{cor}
Clearly as  $z$  approaches unity in this we have (5.9), or if $n$
is increased indefinitely (with  $|z| < 1$) we have essentially
the logarithmic derivative of (4.13). A very interesting corollary
of theorem 5.2 is found from letting  $c_k$  approach zero.  This
results in many simple new identities. It is also interesting to
compare (5.10) with (4.13) to (4.15) inferring a relation between
the Jordan totient function and the Stirling numbers of the second
kind. Another related yet distinct summation formula derived in
the vein of the above analysis is
\begin{thm}   If  $(a_k), (b_k),$   are
each arbitrary sequences chosen so that, together with choice of
$x$, the following functions are defined then
  \begin{align*}
&\sum^n_{k=1} k \ a_k\frac{1-\exp(b_k x)}{1-\exp(b_k x/k)}\\
= &\left(\sum^n_{k=1} a_k\right) +\left(\sum^{[n/2]}_{k=1}
{{\sum_2}_2 a_{2k} \exp((b_{2k} j_1})\frac{x}{2} )\right) +
\left(\sum^{[n/3]}_{k=1} {{\sum_2}_3 a_{3k} \exp((b_{3k}
j_1})\frac{x}{3} )\right)\\
+&\left(\sum^{[n/4]}_{k=1} {{\sum_2}_4 a_{4k} \exp((b_{4k}
j_1})\frac{x}{4} )\right) +\dots + \left(\sum^{[n/n]}_{k=1}
{{\sum_2}_n a_{nk} \exp((b_{nk} j_1})\frac{x}{n} )\right).
  \end{align*}
\end{thm}
This result with $a_k=z^k/k, \, b_k=k(logx)/x, \,
n\longrightarrow\infty$  is
\begin{cor}  If $|z|, |xz|$ are both $< 1$ and $_m \varphi_{k}$  is the number of solutions
in integers of the conditions
 \begin{equation}
(a,m,k) =1, \; 0\leq a \leq k, \; 0\leq m \leq k, \; a+m\neq 0,
 \end{equation}
for fixed  $m$  and  $k$,  then
\begin{equation}
\prod^\infty_{k=1} \left(\frac{1}{1-x^mz^k}\right)^{_m
\varphi_{k}/k} =\exp \left\{ \frac{1}{1-x} \left(\frac{z}{1-z}
-\frac{xz}{1-xz}\right)\right\}.
\end{equation}
\end{cor}
This identity and similarly derived products in higher dimensions
are non-trivial limiting cases of some of the results already
given in Campbell~\cite{gC1998} to~\cite{gC1998a}.  That is, for
example, (5.12) is a limiting case of an identity like (5.3).
    We may write down the generalized version of theorems 5.1 and 5.2 as follows:
\begin{thm}
If  $(a_k), (_1 b_k), (_2 b_k),\ldots (_h b_k)$ are each arbitrary
sequences of functions chosen so that, together with choice of
$x$,  the following series of functions are defined then
\begin{equation} \notag
\sum^n_{k=1}\left(a_k \prod^h_{\lambda =1}\frac{1-\exp(_\lambda
b_k x)}
{1-\exp(_\lambda b_k x/k)}\right)\\
  \end{equation}
  \begin{equation} \notag
= \left(\sum^n_{k=1} a_k\right) +\sum^{[n/2]}_{k=1} {\sum_{h}}_2
a_{2k} \exp((_1 b_{2k} j_1 +\dots +_h b_{2k} j_h)\frac{x}{2})\\
  \end{equation}
  \begin{equation} \notag
 +\sum^{[n/3]}_{k=1} {\sum_{h}}_3 a_{3k} \exp((_1 b_{3k}
j_1+\dots +_h b_{3k} j_h)\frac{x}{3})\\
  \end{equation}
  \begin{equation}
 +\sum^{[n/4]}_{k=1} {\sum_{h}}_4 a_{4k} \exp((_1 b_{4k}
j_1+\dots +_h b_{4k} j_h)\frac{x}{4})\\
  \end{equation}
  \begin{equation} \notag
\vdots\\
  \end{equation}
  \begin{equation} \notag
 +\sum^{[n/n]}_{k=1} {\sum_{h}}_n a_{nk} \exp((_1 b_{nk}
j_1+\dots +_h b_{nk} j_h)\frac{x}{n}).
  \end{equation}\end{thm}
We may now write down the general derivative with respect to  $x$
of both sides of this theorem as follows
\begin{cor}  If $m$ is a
positive integer; and $(a_k), (_1 b_k), (_2 b_k),\ldots (_h b_k),$
are each arbitrary sequences of functions chosen so that, together
with choice of $x$, the following series of functions are all
defined then
    \begin{equation} \notag
\sum^n_{k=1} (a_k \, _h P_m(k))\\
  \end{equation}
  \begin{equation} \notag
  = \sum^{[n/2]}_{k=1} a_{2k}{\sum_{h}}_2 ((_1 b_{2k} j_1+\dots +_h
b_{2k} j_h)/2)^m\\
  \end{equation}
  \begin{equation} \notag
  +\sum^{[n/3]}_{k=1} a_{3k} {\sum_{h}}_3 ((_1 b_{3k} J_1+\dots +_h
b_{3k} j_h)/3)^m\\
  \end{equation}
  \begin{equation}
  +\sum^{[n/4]}_{k=1} a_{4k} {\sum_{h}}_4 ((_1
b_{4k} J_1+\dots +_h b_{4k} j_h)/4)^m\\
  \end{equation}
  \begin{equation} \notag
  \vdots\\
  \end{equation}
  \begin{equation} \notag
  +\sum^{[n/n]}_{k=1} a_{nk} {\sum_{h}}_n ((_1 b_{nk} J_1+\dots +_h
b_{nk} j_h)/n)^m,
\end{equation}
where  $_h P_m (k)$ is the coefficient of  $x^m$  in
 \begin{equation}
 \prod^h_{\lambda=1} \left\{ \sum^\infty_{\mu=1} T_\mu (_\lambda
 b_k x)^{\mu-1}\right\} \text{ where } T_\mu=-\sum^\mu_{\alpha=1}
 \left(\begin{matrix} \mu\\\alpha \end{matrix}\right)
 \frac{B_\alpha}{k^{\alpha-1}}
 \end{equation}
with  $B_\alpha$  the Bernoulli numbers.
\end{cor}
Theorem 5.5 results from theorem 5.10 if  $x$  approaches zero,
whilst corollary 5.11 is the result of equating coefficients of
powers of  $x$  from expanding the two sides of theorem 5.10 as
power series in  $x$.  There is no problem about this as both
sides of (5.14) are finite sums.  The formally complex nature of
$_h P_m (k)$ means that the simplest cases of corollary 5.11 are
those in which  $h = 2$  and  $m = 1, 2, 3, \ldots$  successively.
Clearly,

\begin{equation}
{}_2 P_1(k)=2 T_1 T_2 \ _1 b_k \ _2 b_k,
 \end{equation}
\begin{equation}\notag
{}_2P_2(k)=T_1T_3 (_1b_k)(_2b_k)^3 +T_2 T_2(_1b_k)^2 (_2b_k)^2+T_3
T_1(_1b_k)^3 (_2b_k)
 \end{equation}
\begin{equation}
=T_1T_3(_1b_k)(_2b_k)((_1b_k)^2+(_2b_k)^2)+(T_2)^2(_1b_k)^2(_2b_k)^2,
\end{equation}

so applying (5.16) and then (5.17) to corollary 5.11 gives us
successively the two following corollaries.
\begin{cor}If $(a_k), (_1 b_k), (_2 b_k)$
are each arbitrary sequences of functions chosen so that the
following series of functions are all defined then for  $n$
positive integers greater than  $1$,

\begin{equation} \notag
 \frac{1}{3} \sum^n_{k=1} \left(\frac{1}{k} a_k
{}_1b_k\right) =\sum^{[n/2]}_{k=1} \frac{1}{2} a_{2k} {\sum_{2}}_2
(_1b_{2k}
j_1+_2b_{2k} j_2)\\
 \end{equation}
\begin{equation} \notag
+\sum^{[n/3]}_{k=1} \frac{1}{3} a_{3k} {\sum_{2}}_3 (_1b_{3k}
j_1+_2b_{3k} j_2)\\
 \end{equation}
\begin{equation}
+\sum^{[n/4]}_{k=1} \frac{1}{4} a_{4k} {\sum_{2}}_4 (_1b_{4k}
j_1+_2b_{4k} j_2)\\
 \end{equation}
\begin{equation} \notag
\vdots\\
 \end{equation}
\begin{equation} \notag
+\sum^{[n/n]}_{k=1} \frac{1}{n} a_{nk} {\sum_{2}}_n (_1b_{nk}
j_1+_2b_{nk} j_2).
 \end{equation}
\end{cor}

\begin{cor}  For the same conditions as corollary 5.12,
 \begin{align*}
& \frac{1}{2} \sum^n_{k=1} a_k \left\{ \left(\frac{1}{6}
k^2-\frac{1}{4} k+\frac{1}{12}\right)(_1b_k^2+_2b_k^2)
\frac{1}{4} (k-1)^2 {}_1b_k {} \, _2b_k\right\}\\
&=\sum^{[n/2]}_{k=1} \frac{1}{2} a_{2k} {\sum_{2}}_2 (_1b_{2k}
j_1+_2b_{2k} j_2)\\
&+\sum^{[n/3]}_{k=1} \frac{1}{3} a_{3k} {\sum_{2}}_3 (_1b_{3k}
j_1+_2b_{3k} j_2)\\
&+\sum^{[n/4]}_{k=1} \frac{1}{4} a_{4k} {\sum_{2}}_4 (_1b_{4k}
j_1+_2b_{4k} j_2)\\
&\vdots\\
&+\sum^{[n/n]}_{k=1} \frac{1}{n} a_{nk} {\sum_{2}}_n (_1b_{nk}
j_1+_2b_{nk} j_2).
 \end{align*}
\end{cor}

It is clear that cases of our function $\phi_t(m;k)$ from
definition 4.2 are a natural occurrence in the right sides of
corollaries 5.12 and 5.13
 if $_1b_k=1$ and $_2b_k=1$. When this is the case we
have the two corollaries,

\begin{cor}  If $(a_k)$ is an arbitrary sequence of functions chosen so that
the following series of functions are all defined, and
$\phi_m(n)$ is the $m$th  power of the sum of the non-negative
integers  $a, b$, less than  $n$  such that  $(a, b, n) = 1$  and
$a + b \neq 0$, then for positive integers  $n > 1$,

 \begin{equation}
\sum^n_{k=1} a_k k(k-1)=\frac{1}{2} \phi_1(2) \sum^{[n/2]}_{k=1}
a_{2k} +\frac{1}{3}\phi_1(3) \sum^{[n/3]}_{k=1} a_{3k} +\dots +
\frac{1}{n}\phi_1(n) \sum^{[n/n]}_{k=1} a_{nk},
 \end{equation}

\begin{equation}
\sum^n_{k=1} a_k (\frac{7}{12} k^2-k+\frac{5}{12})=\frac{1}{2}
\phi_2(2) \sum^{[n/2]}_{k=1} a_{2k} +\frac{1}{3}\phi_2(3)
\sum^{[n/3]}_{k=1} a_{3k} +\dots + \frac{1}{n}\phi_2(n)
\sum^{[n/n]}_{k=1} a_{nk},
 \end{equation}
\end{cor}

We next state some examples using this corollary.  The cases are
fairly obvious given the previous analysis,
    \begin{cor}

 \begin{equation}
\sum^n_{k=1} k(k-1)=\frac{1}{2}\phi_1(2)
[n/2]+\frac{1}{3}\phi_1(3)[n/3]+\dots +\frac{1}{n}\phi_1(n)[n/n],
 \end{equation}

\begin{equation}
\sum^n_{k=1} (\frac{7}{12}k^2-k+\frac{5}{12})=\frac{1}{2}\phi_2(2)
[n/2]+\frac{1}{3}\phi_2(3)[n/3]+\dots +\frac{1}{n}\phi_2(n)[n/n],
 \end{equation}

\begin{multline}
\sum^n_{k=1}  k(k-1)\\ =\frac{1}{2}\phi_1(2) [n/2]
(1+[n/2])+\frac{1}{3}\phi_1(3)[n/3](1+[n/3])+\dots
+\frac{1}{n}\phi_1(n)[n/n](1+[n/n]),
 \end{multline}

 \begin{multline}
\sum^n_{k=1}k\left(\frac{7}{12}k^2-k+\frac{5}{12}\right)\\
=\frac{1}{2}\phi_2(2)
[n/2](1+[n/2])+\frac{1}{3}\phi_2(3)[n/3](1+[n/3])+\dots
+\frac{1}{n}\phi_2(n)[n/n](1+[n/n]),
 \end{multline}
\end{cor}

These results are akin to those in Campbell~\cite{gC1994a} and
Ramanujan~\cite{sR1962}.  Furthermore, if we permit  $n$  to
increase indefinitely and let  $a_k=k^{-(s+2)}$   in  corollaries
5.13 and 5.14 we have the Dirichlet generating functions given by
\begin{cor}  If  $\Re s >1$   then

\begin{equation}
\frac{\zeta (s)-\zeta(s+1)}{\zeta(s+2)}=\sum^\infty_{k=2}
\frac{\phi_1(k)}{k^{(s+3)}},
\end{equation}

\begin{equation}
\frac{7\zeta (s)-12
\zeta(s+1)+5\zeta(s+2)}{\zeta(s+3)}=\sum^\infty_{k=2}
\frac{\phi_2(k)}{k^{(s+3)}},
\end{equation}
\end{cor}

Likewise, if  $n$  too increases indefinitely and $a_k=z^k/k$   in
corollaries 5.13 and 5.14 we have the infinite products,
    \begin{cor}  If   $|z|<1$  then

\begin{equation}
\prod^\infty_{k=2} (1-z^k)^{-\phi_1(k)/k^2}
=\exp\left\{\frac{z}{(1-z)^2}\right\},
\end{equation}

 \begin{equation}
\prod^\infty_{k=2} (1-z^k)^{-\phi_2(k)/k^2}
=\left(\frac{1}{1-z}\right)^{\frac{5}{12}}
\exp\left\{\frac{z(12z-5)}{12(1-z)^2}\right\}.
\end{equation}
\end{cor}

    Recalling (4.11) and its ensuing paragraph, it was clear that

 \begin{equation}
\frac{\zeta (s-m)}{\zeta(s)}=1+\sum^\infty_{k=2}
\frac{J_m(k)}{k^s}, \qquad \mathcal{R} s>1.
\end{equation}
This sum when compared to corollary 5.6 shows that
    \begin{cor}  For positive integers  $k$,

 \begin{equation}
\phi_1(k) =J_2(k)-J_1(k),
 \end{equation}

\begin{equation}
\phi_2(k)=\frac{7}{12}J_3(k) -J_2(k)+\frac{5}{12}J_1(k).
\end{equation}
\end{cor}

\section{Application of Jacobi theta series to the generalized
summations.}

In Campbell~\cite{gC1993}, the Jacobi theta function was applied
to the finite left hand side form of the summation formula (3.2).
We use now the well known terminology for the theta function

 \begin{equation}
\Theta_1(z,q)=\Theta_1(z)=2q^{\frac{1}{2}} \sum^\infty_{k=1}
(-1)^k q^{k(k+1)}\sin (2k+1)z,
 \end{equation}

\begin{equation}
\sum^\infty_{k=1} \frac{1}{k} \frac{q^{2k}}{1-q^{2k}} \sin 2k
\alpha \sin 2k\alpha =\frac{1}{4}
\log\left(\frac{\Theta_1((\alpha+\beta),q)\sin(\alpha
-\beta)}{\Theta_1((\alpha -\beta),q)\sin (\alpha+\beta)}\right).
\end{equation}

Application of (6.2) to (3.2) gives us the result (see Campbell
~\cite[section 4]{gC1993})
    \begin{thm}  If  $|q|<1$ and $x^{1/k}\neq 1$   then

 \begin{multline}
\exp\left(4\sum^\infty_{k=1} \frac{1}{k}
\frac{q^{2k}}{1-q^{2k}}\frac{1-x}{1-x^{1/k}}\sin 2k\alpha \sin 2k
\alpha\right)\\
=\frac{\Theta_1(\alpha+\beta) \sin
(\alpha-\beta)}{\Theta_1(\alpha-\beta) \sin (\alpha+\beta)}
\prod^\infty_{k=1}\left(\frac{\Theta_1((\alpha+\beta)
k,q^{k})\sin(\alpha -\beta)}{\Theta_1((\alpha -\beta)k,q^{k})\sin
(\alpha+\beta)}\right)^{f_k(x)/k}
 \end{multline}

where for positive integers $j$ less than $k$,
\begin{equation}\notag f_k(x)= \sum_{(j,k)=1} x^{j/k}=
{\sum_{1}}_k x^{j_1/k}.\end{equation}
\end{thm}

Particular cases of this statement include some of the known
arithmetical functions, as highlighted in the following examples.
To begin with, for positive integers  $n$  with $x=exp(2\pi in)$
we have
    \begin{cor}  (see Campbell [10])

\begin{multline}
\exp\left(4\sum^\infty_{k\mid n} \frac{q^{2k}}{1-q^{2k}} \sin
2k\alpha \sin 2k
\alpha\right)\\
=\frac{\Theta_1(\alpha+\beta) \sin
(\alpha-\beta)}{\Theta_1(\alpha-\beta) \sin (\alpha+\beta)}
\prod^\infty_{k=1}\left(\frac{\Theta_1((\alpha+\beta)
k,q^{k})\sin(\alpha -\beta)}{\Theta_1((\alpha -\beta)k,q^{k})\sin
(\alpha+\beta)}\right)^{c_k(n)/k},
 \end{multline}

where  $c_k(n)$   is Ramanujan's trigonometrical function.
\end{cor}
    If in (6.3) we allow  $x$  to approach unity, we have the result,
    \begin{cor}

\begin{multline}
\exp\left(4\sum^\infty_{k=1} \frac{q^{2k}}{1-q^{2k}} \sin 2k\alpha
\sin 2k
\alpha\right)\\
=\frac{\Theta_1(\alpha+\beta) \sin
(\alpha-\beta)}{\Theta_1(\alpha-\beta) \sin (\alpha+\beta)}
\prod^\infty_{k=1}\left(\frac{\Theta_1((\alpha+\beta)
k,q^{k})\sin(\alpha -\beta)}{\Theta_1((\alpha -\beta)k,q^{k})\sin
(\alpha+\beta)}\right)^{\varphi(n)/k},
 \end{multline}

where  $\varphi(n)$  is the Euler totient function.
\end{cor}

If we follow an analogous line of reasoning to that which led us
to (6.3) but from starting with lemma 3.1, we obtain the new
result,
    \begin{thm}  If  $|q|<1$ and all of $x_1^{1/k}, x_2^{1/k},\ldots,x_m^{1/k}\neq1$   then

\begin{multline}
\exp\left(4\sum^\infty_{k=1} \frac{1}{k} \frac{q^{2k}}{1-q^{2k}}
\left(\frac{1-x_1}{1-x_1^{1/k}} \frac{1-x_2}{1-x_2^{1/k}} \dots
\frac{1-x_m}{1-x_m^{1/k}}\right)\sin 2k\alpha \sin 2k
\alpha\right)\\
=\frac{\Theta_1(\alpha+\beta) \sin
(\alpha-\beta)}{\Theta_1(\alpha-\beta) \sin (\alpha+\beta)}
\prod^\infty_{k=1}\left(\frac{\Theta_1((\alpha+\beta)
k,q^{k})\sin(\alpha -\beta)}{\Theta_1((\alpha -\beta)k,q^{k})\sin
(\alpha+\beta)}\right)^{_m f_k(x)/k}
 \end{multline}

where  \begin{equation}_m f_k(x)= {\sum_m}_k \left( {x_1}^{j_1}
{x_2}^{j_2} \cdots {x_m}^{j_m} \right)^{1/k}.\end{equation}
\end{thm}

We state two relevant corollaries with respect to the substance of
the current work. Firstly, to obtain the corresponding result for
the Jordan totient function, we need only allow all of the $x$'s
to approach unity, bearing in mind our work in section~4. We
therefore have
    \begin{cor}

 \begin{multline}
\exp\left(4\sum^\infty_{k=1} k^{m-1}\frac{q^{2k}}{1-q^{2k}} \sin
2k\alpha \sin 2k \alpha\right)\\
=\frac{\Theta_1(\alpha+\beta) \sin
(\alpha-\beta)}{\Theta_1(\alpha-\beta) \sin (\alpha+\beta)}
\prod^\infty_{k=1}\left(\frac{\Theta_1((\alpha+\beta)
k,q^{k})\sin(\alpha -\beta)}{\Theta_1((\alpha -\beta)k,q^{k})\sin
(\alpha+\beta)}\right)^{J_m (k)/k},
 \end{multline}

where  $J_m (k)$  is the Jordan totient function.
\end{cor}

Secondly, to obtain the corresponding result for our new function,
$c_k(n_1, n_2, \dots, n_m)$ we need only allow all of the $x_j$'s
to approach $exp(2\pi i n_j)$, bearing in mind our work in
section~1. We therefore have
    \begin{cor}

 \begin{multline}
\exp\left(4\sum_{k|(n_1, n_2, \dots, n_m)} \frac{q^{2k}}{1-q^{2k}}
\sin
2k\alpha \sin 2k \alpha\right)\\
=\frac{\Theta_1(\alpha+\beta) \sin
(\alpha-\beta)}{\Theta_1(\alpha-\beta) \sin (\alpha+\beta)}
\prod^\infty_{k=1}\left(\frac{\Theta_1((\alpha+\beta)
k,q^{k})\sin(\alpha -\beta)}{\Theta_1((\alpha -\beta)k,q^{k})\sin
(\alpha+\beta)}\right)^{c_k(n_1, n_2, \dots, n_m)/k}.
 \end{multline}

where  $c_k(n_1, n_2, \dots, n_m)$  is our new generalized
Ramanujan totient function.
\end{cor}

We have not attempted to use the full generality of corollary 3.2,
and it is clear that we may return to this topic in future papers.
The return to such investigations will no doubt be hastened if it
turns out that the identities given in this paper can be shown to
enumerate tilings in the sense described by
Baake~et~al~\cite{mB1994} and his colleagues pertaining to
quasicrystals.




\end{document}